\documentclass{amsart}
\usepackage{amssymb}
\usepackage{graphicx}

\usepackage{amsmath,bm}
\usepackage{amsthm}
\usepackage{amsfonts}
\usepackage{latexsym}
\usepackage{color}
\usepackage{stmaryrd}
\usepackage[normalem]{ulem}
\usepackage{setspace}

\newtheorem{theorem}{Theorem}[section]

\theoremstyle{definition}

\theoremstyle{remark}

\numberwithin{equation}{section}

\newcommand{\R}{\mathbb{R}}

\begin{document}

\title[Why condensation by compression cannot occur]{Why condensation by compression in pure water vapor cannot occur in an approach based on Euler equations}


\author{Maren Hantke}
\address{Institute for Analysis and  Numerics,
Otto-von-Guericke-University Magdeburg, PSF 4120, D--39016 Magdeburg,
Germany.}
\curraddr{Maren Hantke, Institute for Applied Mathematics, Gottfried Wilhelm Leibniz University Hannover, Welfengarten 1, D--30167 Hannover, Germany.}
\email{maren.hantke@ovgu.de}

\author{Ferdinand Thein}
\address{Institute for Analysis and  Numerics,
Otto-von-Guericke-University Magdeburg, PSF 4120, D--39016 Magdeburg,
Germany.}
\email{ferdinand.thein@ovgu.de}
\thanks{The second author was supported by \textit{Landesgraduiertenstipendium Sachsen-Anhalt} and 
is supported by the DFG grant HA 6471/2-1. The authors thankfully acknowledge the support.}


\date{}


\begin{abstract}
Phase transitions are in the focus of the modeling of multiphase flows. A large number of models is available to describe such processes. 
 We consider several different two phase models that are based on the Euler equations of compressible fluid flows and which take into 
 account phase transitions between a liquid phase and its vapor. Especially we consider the flow 
 of liquid water and water vapor. We give a mathematical proof, that all these models are not able to describe the process 
 of condensation by compression. This behavior is in agreement with observations in experiments, that simulate adiabatic flows, and shows that 
the Euler equations give a fairly good description of the process. The mathematical proof is valid for the official standard {\em IAPWS-IF97} 
for water and for any other 
 good equation of state. Also the opposite case of expanding the liquid phase will be discussed.
\end{abstract}

\maketitle

\section{Introduction} 
Compressible two and multi phase flows occur in nature as well as in numerous industrial applications. In many cases 
phase transitions between the fluids are of major importance. Examples are the formation of clouds or the phenomenon of 
cavitation that for instance appears in liquid water flows around ship propellers.

The modeling of such processes is a challenge. The description of the interfaces between the fluids as well as their 
interaction is of high complexity. Therefore in the spotlight of the methods is the treatment of the interfaces. 
Many numerical simulations are based on the Euler equations of compressible fluid flow. We 
will direct our attention to two phase models of this type that take into account mass transfer between the fluids.

Very famous are the models of Baer-Nunziato type. Here both phases are described by their own set of Euler equations. An 
additional equation for the volume fractions of the phases is considered, see Section \ref{s:moddis}. The original model of Baer and 
Nunziato \cite{Baer-Nunziato:1986}, that does not include the effect of mass exchange between the phases, was modified 
by Saurel and Abgrall in \cite{Saurel-Abgrall:1999} by introducing relaxation terms for pressures and velocities of the 
phases. Later in \cite{Saurel-Petitpas-Abgrall:2008} a similar idea allowed the description of phase transition by 
using relaxation terms for the temperatures and chemical potentials. This idea was picked up for instance by 
Petitpas et al.\ in \cite{Petitpas-usw:2009} or by Zein et al. in \cite{Zein-Hantke-Warnecke:2010}.

Another type of modeling of two phase flows is to use only one set of Euler equations. 
Each phase has its own equation of state. Phase transitions can be described by a further equation that is called kinetic 
relation. See for instance the well reputed article of Abeyaratne and Knowles \cite{Abeyaratne} that deals with solid-solid 
interfaces, the papers of Merkle \cite{Merkle} or Hantke et al. \cite{Hantke} on the isothermal Euler equations.

Finally we want to refer to a recent paper of Dumbser et al.\ \cite{Dumbser}. In their work also only one set of 
equations is used. Phase transitions take place only in thermal equilibrium, no kinetic relation is used. 
Surprisingly this type of modeling is closely related to the Baer-Nunziato type modeling mentioned above including 
relaxation terms to describe mass transfer. We come back 
to this and the abovementioned models later in Section \ref{s:moddis}.

One can find an extensive literature on cavitating flows, but the opposite question of the creation of a liquid phase by a strong compression  
is discussed only in rare cases. In the following we consider pure water vapor, that will be highly compressed. This can be realized 
by a steam filled tube with a flexible piston, see Figure \ref{comkolb}. 
If there 
is no heat exchange with the neighbourhood of the tube, the process is nearly an adiabatic flow. Therefore it can be fairly described
by the compressible Euler equations. One may expect, 
that it is possible to compress the vapor phase such that the vapor will condensate. This means that a liquid phase is created. 
In fact, it turns out this is impossible in a non-isothermal approach based on Euler equations, which is in agreement with observations from experiments. 
The main focus of our work is to give a {\bf mathematical proof for this phenomenon}. 

In the case of expanding a liquid under the same boundary conditions, see Figure \ref{expkolb} 
the situation is more complex. Nevertheless, also for the cavitation case 
we can prove some theoretical results. For detailed discussions of cavitation models we refer to Iben \cite{I2} and \cite{I1}.

The paper is organized as follows. In Section \ref{s:eulergl} we introduce the compressible Euler equations and 
briefly recall some well-known analytical results on the Riemann problem. Afterwards in Section \ref{s:eos} we 
give some equations of state to close the system and provide some physical background. Section \ref{s:nucl} deals 
with compressed water vapor. First in Section \ref{s:idea} we explain the idea for the proof of our statement. 
In Section \ref{s:spec} we show for a special choice 
of equations of state for the phases, that condensation by compression 
cannot occur. 
After that in Section \ref{s:moddis} we show, that this idea is also applicable to Baer-Nunziato type models 
with relaxation terms as in Saurel et al.\ \cite{Saurel-Petitpas-Abgrall:2008}, Zein et al.\ 
\cite{Zein-Hantke-Warnecke:2010}. Thereafter, in Section \ref{s:exteos} we generalize the proof to the {\em real equation of state for water}. 
Finally we consider the opposite case of cavitation by expansion in Section \ref{s:cav}. 
We end up with some closing remarks in 
Section \ref{s:clos}.


\section{The Euler model}\label{s:eulergl}
The compressible Euler equations in one space dimension are given by the following 
system
\begin{eqnarray}
\label{e:mb}
\frac{\partial }{\partial t}\rho+\frac{\partial }{\partial x}\rho u&=&0\,,\\
\label{e:ib}
\frac{\partial }{\partial t}\rho u+\frac{\partial }{\partial x}(\rho u^2+p)&=&0\,,\\
\label{e:eb}
\frac{\partial }{\partial t}\rho(e+\frac{1}{2}u^2)+\frac{\partial }{\partial x}\bigl[\rho(e+\frac{1}{2}u^2)+p\bigr]u&=&0\,,
\end{eqnarray}
where the variables $\rho$, $u$ and $e$ denote the {\em mass density}, the {\em velocity} and the {\em internal energy}, 
resp. The further quantity $p$ describes the {\em pressure}. It is related to the mass density $\rho$ and the internal 
energy $e$ by an {\em equation of state}, see Section \ref{s:eos}. All physical fields depend on time $t\in\R_{\ge0}$ 
and on space $x\in\R$.

Here we consider Riemann problems for the Euler equations, that are given by the above balance equations 
(\ref{e:mb})-(\ref{e:eb}), an equation of state and the corresponding {\em Riemann initial data}
\begin{equation}\label{e:RI}
\rho(0,x)=\left\{
\begin{array}{c}
\rho_-\\\rho_+
\end{array}
\right.\,,
\quad
u(0,x)=\left\{
\begin{array}{c}
u_-\\u_+
\end{array}
\right.\,,
\quad
e(0,x)=\left\{
\begin{array}{ccl}
e_-&\mbox{for}&x<0\\e_+&\mbox{for}&x\ge0\,.
\end{array}
\right.
\end{equation}
This is the simplest choice of initial conditions with piecewise constant data. It is possible and 
conventional to give initial states  for $(\rho, u, p)$ or $(p, u, T)$ instead of initial states for $(\rho, u, e)$.
Also other choices are imaginable.

The Riemann problem is very helpful in the context 
of hyperbolic partial differential systems, because it exhibits all phenomena as shock or rarefaction waves. 
It is a basic problem in the theory of hyperbolic systems. 
In numerics Riemann problems appear in finite volume methods for systems of conservation laws due to the discreteness 
of the grid.

For the Riemann problem for the compressible Euler equations equipped with an appropriate equation of state 
one can construct the exact solution. The solution is selfsimilar. It consists of constant states, that are separated 
by shock and rarefaction waves and a contact discontinuity. Details can be found in several textbooks. For basics 
on conservation laws see the books of Toro \cite{TORO}, Lax \cite{LAX}, 
LeVeque \cite{bLEV}, Smoller \cite{smo}, Kr\"oner \cite{kro}, Dafermos \cite{DAF}, 
and others.


\section{Equation of state}\label{s:eos}

As mentioned in the previous section we need an equation of state to close the system 
(\ref{e:mb})-(\ref{e:eb}). Several commonly accepted equations of state are available like the van der 
Waals equation of state or the Tait's equation. A collective problem is that for any choice of parameters 
all these equations 
at the best only locally give a good approximation of the thermodynamic properties of water vapor or liquid water.

On the other hand, the real equation of state for water to the official standard {\em IAPWS-IF97} based on the standard 
formulation of Wagner et al.\ \cite{Wagner-Kruse:1998}, \cite{Wagnerusw:2000}, \cite{steam}
is too complex for analytical consideration. In the following this equation of state is called {\em real equation of state}.

For the moment we use a modified form of the stiffened gas equation of state, 
see \cite{Saurel-Petitpas-Abgrall:2008}.
It is given by the following relations
\begin{eqnarray}\label{e:eos-e}
e_k(p_k,\rho_k)&=&\frac{p_k+\gamma_k\pi_k}{\rho_k(\gamma_k-1)}+q_k\,,\\
\label{e:eos-T}
T_k(p_k,\rho_k)&=&\frac{p_k+\pi_k}{C_k\rho_k(\gamma_k-1)}\,,\\
\label{e:eos-a}
a_k(p_k,\rho_k)&=&\sqrt{\frac{\gamma_k(p_k+\pi_k)}{\rho_k}}\,,\\
\label{e:eos-s}
s_k(p_k,T_k)&=&C_k\ln\frac{T_k^{\gamma_k}}{(p_k+\pi_k)^{(\gamma_k-1)}}+q_k'\,.
\end{eqnarray}
Here $T$ and $s$ denote the {\em temperature} and the {\em specific entropy} of the fluid. The speed of sound is given by $a$.
The index $k=V,L$ indicates the fluid 
under consideration, vapor or liquid. The parameters $\gamma$, $\pi$, $q$, $q'$ and $C$ will be specified later. 
Note, that for the special choice of $\pi=0$ and $q=0$ the equation of state reduces to the ideal gas law.

For the physical background of the following considerations and more details on 
thermodynamics we refer to the books of M\"uller and M\"uller 
\cite{bMUE} and M\"uller \cite{bMUE1}.

The {\em specific Gibbs free energy} of the phases is given by
\begin{equation}\label{e:g}
g_k=e_k+\frac{p_k}{\rho_k}-T_ks_k\,.
\end{equation}
If the vapor and the liquid phase are in thermodynamic equilibrium the Gibbs free energies equal each other, 
this means that 
\begin{equation}\label{eqcon}
g_V=g_L\,.
\end{equation}
Using the relations (\ref{e:eos-e})-(\ref{e:eos-s}) the Gibbs free energy of each phase can be expressed as 
a function of the temperature and the pressure 
\[
g_k=g_k(p_k,T_k)\,.
\]
Then the equilibrium condition (\ref{eqcon}) gives a direct relation between the temperature and the corresponding 
equilibrium pressure or {\em saturation pressure}
\begin{equation}\label{psat}
p_{sat}=p_{sat}(T)\,.
\end{equation}
This relation gives the {\em saturation line} in the $(T,p)$-phase space. Sometimes it is useful to inversly express 
the temperature as a function of the pressure
\begin{equation}\label{e:Tsat}
 T_{sat}=T_{sat}(p)\,.
\end{equation}
For admissible pressures we obtain the corresponding {\em saturation temperature}.

\begin{figure}[h!]
\includegraphics[scale=0.6]{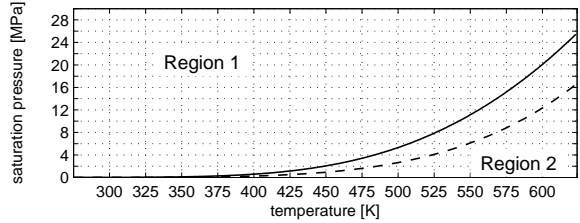}
\caption{{\footnotesize{Solid line: saturation curve $p_{sat(T)}$ for parameters given in Table \ref{tab:para}, 
dashed line: real saturation curve}}}
\label{fig:sat}
\end{figure}

For the moment we use the same parameters as Saurel et al.\ \cite{Saurel-Petitpas-Abgrall:2008}. These parameters 
are given in Table \ref{tab:para}.
\begin{table}[h!]
\caption{Parameters for water vapor and liquid water, \cite{Saurel-Petitpas-Abgrall:2008}}
\label{tab:para}
\begin{tabular}{ccccccc}
\hline
$k$&$\gamma$&$\pi$ [Pa]&$C$ [J/kg/K]&$q$ [J/kg]&$q'$ [J/kg/K]&\\
vapor&1.43&0&1040&2030000&-23000\\
liquid&2.35&$10^9$&1816&-1167000&0\\
\hline
\end{tabular}
\end{table}
For this special choice of parameters
we obtain the saturation curve given by the solid line in Figure \ref{fig:sat}.

Here Region 1 belongs to the liquid water phase, whereas Region 2 belongs 
to the water vapor. The dashed line marks the real saturation line.
Obviously the precise shape of the saturation line directly depends on the choice for the equations of state and 
the parameters therein.

Note that for a thermodynamic consistent equation of state the following relations hold
\begin{equation}\label{e:gibbsrel}
\frac{\partial g(p,T)}{\partial p}=\frac{1}{\rho} \qquad \mbox{and} \qquad \frac{\partial g(p,T)}{\partial T}=-s\,.
\end{equation}
These relations will be used for the proof of our statement. 


\section{Condensation by compression}\label{s:nucl}
\subsection{Wave curve in the $(p,T)$-phase space}\label{s:idea}
\begin{figure}[h!]
\includegraphics[scale=0.4]{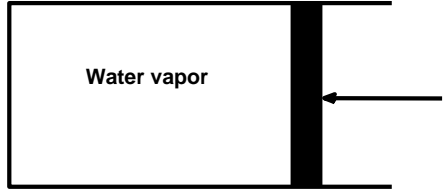}
\caption{{\footnotesize{Compression of water vapor}}}
\label{comkolb}
\end{figure}

At first we consider the case of the compression of water vapor. This case can be simulated by a steam filled tube 
equipped with a flexible piston, which is highly sped up to compress the vapor phase, see Figure \ref{comkolb}.
The compression of water vapor will lead to an increase of the pressure and the density of the vapor phase. 
From the theory of the Euler equations we know, that a shock wave will propagate through the vapor phase. 

Assume that the state ahead of the shock is given by $(\hat \rho, \hat u, \hat p)$. The state behind the shock is 
denoted by 
$(\rho_*,u_*,p_*)$. Then we have the following relationship for the density and the pressure
\begin{equation}\label{e:p_rho}
 \frac{\rho_*}{\hat\rho}=\frac{\left(\frac{p_*}{\hat p}\right)+\left(\frac{\gamma-1}{\gamma+1}\right)
 +\frac{2\gamma\pi}{\hat p(\gamma+1)}}{\left( 
 \frac{\gamma-1}{\gamma+1}\right)\left(\frac{p_*}{\hat p}\right)+1+\frac{2\gamma\pi}{\hat p(\gamma+1)}}\,.
\end{equation}
For details of the derivation of relation (\ref{e:p_rho}) we recommend the book of Toro \cite{TORO}, Section 3.1. 
Note, that this relation holds only for the generalized stiffened gas law (\ref{e:eos-e})-(\ref{e:eos-s}). Using the 
equation of state (\ref{e:eos-T}) we easily obtain an analogous relation for the pressure and the temperature, which 
is given by
\begin{equation}\label{e:p_T}
 \frac{\hat T}{T_*}=\frac{\hat p+\pi}{p_*+\pi}\cdot\frac{p_*(\gamma+1)+\hat p(\gamma-1)+2\gamma\pi}{
 \hat p(\gamma+1)+p_*(\gamma-1)+2\gamma\pi\frac{\hat p}{p_*}}\,.
\end{equation}
From Equation (\ref{e:p_T}) we find the wave curve $T_*(p_*;\hat p,\hat T)$ in the $(p,T)$-phase space, that denotes all states 
$(p_*,T_*)$ that can be connected to the initial state $(\hat p, \hat T)$ by a shock wave.

\begin{figure}[h!]
\includegraphics[scale=0.6]{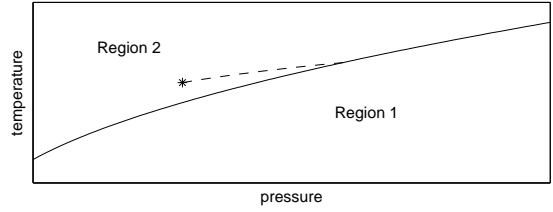}
\caption{{\footnotesize{Solid line: saturation curve $T_{sat}(p)$, dashed line: wave curve, star: initial state}}}
\label{fig:nucl}
\end{figure}

Assume that the vapor phase with initial pressure and temperature $(\hat p,\hat T)$ is compressed 
sufficiently such that a liquid phase is created, 
then the corresponding wave curve $T_*(p_*;\hat p,\hat T)$ must have an intersection point with the saturation line, see 
the sketch in Figure \ref{fig:nucl}.
We want to prove that this is impossible. For the proof we assume the existence of the intersection point and 
derive a contradiction.


\subsection{Proof of the statement for a particular equation of state}\label{s:spec}
Assume, the intersection point of the wave curve and the saturation line is given by $(p_*,T_*)$. Then the function 
\begin{equation}\label{e:Trel}
\hat T(\hat p;p_*,T_*)=T_*\cdot\frac{\hat p}{p_*}\cdot\frac{p_*(\gamma_V+1)+\hat p(\gamma_V-1)}{
 \hat p(\gamma_V+1)+p_*(\gamma_V-1)}
\end{equation}
denotes all admissible initial states $(\hat p,\hat T)$. Here we already used the fact that for the 
vapor phase we have $\pi=0$. 
Let $'$ denote the derivative of the temperature functions with respect to the pressure.
At the intersection point the following relation must hold
\begin{equation}\label{e:intrel}
 \hat T'(p_*)\le T_{sat}'(p_*)\,.
\end{equation}
By a simple calculation we find
\begin{equation}\label{e:abl1}
 \hat T'(p_*;p_*,T_*)=\frac{T_*}{p_*}\cdot \frac{\gamma_V-1}{\gamma_V}\,.
\end{equation}
To find $T_{sat}'(p_*)$ we start with the equilibrium condition (\ref{eqcon}) and express the Gibbs free energies of the 
phases as functions of $p$ and $T$. We obtain
\begin{equation}\label{e:gpT}
 g_k(p,T)=C_kT\gamma_k+q_k-C_kT\ln \frac{T^{\gamma_k}}{(p+\pi_k)^{\gamma_k-1}}-Tq_k'
\end{equation}
with $k=L,V$. From 
\begin{equation}\label{e:f}
f(p,T_{sat}(p))=g_V(p,T)-g_L(p,T)=0
\end{equation} 
we derive by using the implicit function theorem
\begin{eqnarray*}
 \frac{\partial f}{\partial p}&=&C_V(\gamma_V-1)\frac{T_{sat}}{p}-C_L(\gamma_L-1)\frac{T_{sat}}{p+\pi_L}\,,\\
 \frac{\partial f}{\partial T_{sat}}&=&-C_V\ln \frac{T_{sat}^{\gamma_V}}{(p)^{\gamma_V-1}}-q_V'+C_L
 \ln \frac{T_{sat}^{\gamma_L}}{(p+\pi_L)^{\gamma_L-1}}\\
 &=&C_L\gamma_L-C_V\gamma_V+\frac{q_L-q_V}{T_{sat}}
\end{eqnarray*}
and finally
\begin{equation}\label{e:abl2}
 T_{sat}'(p_*;T_*)=T_*\cdot \frac{\frac{C_V(\gamma_V-1)}{p_*}-\frac{C_L(\gamma_L-1)}{p_*+\pi_L}}{C_V\gamma_V 
 -C_L\gamma_L+\frac{q_V-q_L}{T_*}}<
 T_*\cdot \frac{\frac{C_V(\gamma_V-1)}{p_*}}{C_V\gamma_V }\,.
\end{equation}
This is clear because
\[
 -\frac{C_L(\gamma_L-1)}{p_*+\pi_L}<0 \qquad \mbox{and} \qquad
 -C_L\gamma_L+\frac{q_V-q_L}{T_*}>0\,.
\]

Obviously (\ref{e:abl2}) implies
\begin{equation}\label{e:abl3}
 T_{sat}'(p_*;T_*)<\frac{T_*}{p_*}\cdot \frac{\gamma_V-1}{\gamma_V}=\hat T'(p_*;p_*,T_*)\,.
\end{equation}
This is a contradiction, see (\ref{e:intrel}). Accordingly we have
\begin{theorem}\label{n_sp}
 Using the equations of state (\ref{e:eos-e})-(\ref{e:eos-s}) and the parameters given in Table \ref{tab:para} 
 condensation by compression of pure water vapor cannot occur.
\end{theorem}

\subsection{Short discussion of different phase transition models}\label{s:moddis}
In the following paragraph we explain, why the proof of our statement given in Section \ref{s:spec} resp.\
in Section \ref{s:exteos} is applicable for all considered models based on Euler equations.

A pure water vapor phase can be described by a single set of Euler equations. The compression leads to an 
increase of the density and the pressure as already mentioned in the previous section. One may assume that for 
sufficiently strong compression the vapor phase starts to condensate. This means, that a liquid phase is created.

Of course, for any state $(p,T)$ which is in the interior of Region 2 (water vapor), the vapor phase is 
situated in a stable state. In order that phase condensation can happen, there must be a mechanism for phase transition. 
Therefore, it is clear that in the case of condensation the wave curve, see the next to last section, must have an 
intersection point with the saturation line. For this it doesn't matter, whether phase transition is modeled 
by a kinetic relation \cite{Hantke}, \cite{Merkle} or by using an equilibrium assumption as done in \cite{Dumbser}. 

Using a kinetic relation, a nucleation criterion is used, see \cite{Hantke}. 
Here a critical state is reached, in which the vapor phase starts 
to condensate, which implies the intersection point.

In Dumbser et al.\ \cite{Dumbser} phase transition is modeled by an equilibrium assumption. For any given temperature 
$T$ and $p<p_{sat}(T)$ the pair $(T,p)$ describes some vapor state. Analogously $(T,p)$ with $p>p_{sat}(T)$ describes 
the fluid in the liquid state. For $p=p_{sat}(T)$ one may have water vapor or liquid water as well as a mixture 
of both fluids. The fluid at the saturation state is defined by its temperature and the mass fraction or equivalently by its 
pressure and the mass fraction of the vapor/liquid phase. All corresponding states in the $(p,T)$-phase plane 
are located at the saturation line. For more details see Iben et al.\ \cite{Iben}. Nevertheless, for condensation a wave 
curve must have an intersection point with the saturation line.

The Euler equations are only valid for pure fluids or homogeneous mixtures in the
thermodynamic equilibrium. For models that use only one set of Euler equations as discussed before pure fluid are present. 
On the other hand in literatur often models of Baer-Nunziato type are used to describe the situation considered.
The generalized Baer-Nunziato model is given by a 
two phase model using two sets of Euler equations 
\begin{eqnarray*}
\frac{\partial }{\partial t}\alpha_k\rho_k+\frac{\partial }{\partial x}\alpha_k\rho_k u_k&=&\pm \dot m\,,\\
\frac{\partial }{\partial t}\alpha_k\rho_k u_k+\frac{\partial }{\partial x}\alpha_k(\rho_k u_k^2+p_k)&=&\pm P\frac{\partial }{\partial x}\alpha_k\pm M\,,\\
\frac{\partial }{\partial t}\alpha_k\rho_k(e_k+\frac{1}{2}u_k^2)+\frac{\partial }{\partial x}\bigl[\rho_k(e_k+\frac{1}{2}u_k^2)+p_k\bigr]u_k&=&\mp P
\frac{\partial}{\partial t}\alpha_k\pm E\,,
\end{eqnarray*}
$k=1,2$ and a further equation to describe the volume fractions of the phases 
\[
\frac{\partial}{\partial t}\alpha_1+U \frac{\partial }{\partial x}\alpha_1=A
\]
with the same notations as before. Further, $\alpha_k$ denotes the volumefraction of phase $k$. The sources $A,\dot m, M, E$ on the right hand side 
of the equations describe the exchange of mass, momentum and energy. They include relaxation terms for velocity, pressure, 
temperature and Gibbs free energy of the phases, that guarantee, that both phases relax to 
thermodynamic equilibrium. The pressure $P$ und the velocity $U$ have to be defined by some closure law, see \cite{Saurel-Abgrall:1999}, 
\cite{Saurel-Petitpas-Abgrall:2008} or \cite{Zein-Hantke-Warnecke:2010}.

Again we start with 
a pure vapor phase, that will be compressed. For numerical reasons the volume fraction of the pure phase is 
assumed to be $1-\varepsilon$, whereas the volume fraction of the absent phase is assumed to be $\varepsilon$. Typically 
one uses $\varepsilon=10^{-8}$, see for instance Saurel and Abgrall \cite{Saurel-Abgrall:1999}.

Mass transfer is described by the Gibbs free energy relaxation term. Condensation will occur only in the case, 
that the specific Gibbs free energy of the vapor phase is larger than the specific Gibbs free energy of the (artificial) 
liquid phase. This is not the case for any set of initial data that describes a pure water vapor phase. 
Therefore there is no contribution by the relaxation terms as long as $(T,p)$ is in the interior of Region 2 (vapor phase).

In regions of constant volume fractions the system decouples. The solution for each phase can be determined separately. 
This implies that the relations of the single phase Euler equations are also valid for the Baer Nunziato model in the 
case considered. In order that condensation can occur the vapor phase must be compressed in such a manner 
that the specific Gibbs 
free energy of the vapor phase is larger than the specific Gibbs free energy of the (artificial) 
liquid phase. This implies an intersection point of the wave curve with the saturation line.


\section{Extension to the real equation of state for water}\label{s:exteos}

In the previous Section \ref{s:nucl} we have proved, that condensation by compression can not occur for the chosen 
equations of state with parameters given in Table \ref{tab:para}. On the other hand, in Figure \ref{fig:sat} we can see, 
that this choice gives a very bad approximation of the real saturation curve for higher temperatures. We now generalize 
our statement using the results from the last section. We want to show, that for the {\em real equation of state} and for 
{\em good} approximations of the real equation of state condensation by compression cannot occur. The proof uses the same 
arguments as before.

We start with an arbitrary initial state in the vapor region and we consider the corresponding wave curve. 
Assume, $(p_*,T_*)$ is the intersection point of the curve $\hat T(\hat p;p_*,T_*)$ of all admissible initial states 
in the $(p,T)$-phase 
space with the real saturation 
line $T_{sat}(p)$. We compare the derivatives and find the contradiction.

\subsection{Approximation of the real equation of state}\label{approx}
For our purpose it is sufficient to find a good approximation of the real equation of state in a small neighborhood 
of the saturation line.
In the following we show how to find suitable parameters, coming from the intersection point $(p_*,T_*)$. This is an 
improvement of the method of Le Metayer et al. \cite{LeM}, which is a modification of the idea introduced by 
Barberon and Helluy \cite{Bar}. 

\subsubsection{Vapor phase}\label{parv}
For any temperature $T_*$ the corresponding (real) saturation pressure is known by the real formulas given by Wagner 
\cite{Wagnerusw:2000}, \cite{steam}. The same is true for the corresponding vapor density $\rho_{V*}$, 
the speed of sound $a_{V*}$, the entropy $s_{V*}$ and the internal energy $e_{V*}$. For simplicity we choose $\pi_V=0$. 
Then from (\ref{e:eos-a}) we directly obtain $\gamma_V$. Next we calculate $q_V$ from (\ref{e:eos-e}), $C_V$ from 
(\ref{e:eos-T}) and $q_V'$ from (\ref{e:eos-s}).
\begin{figure}[h!]
\includegraphics[scale=0.6]{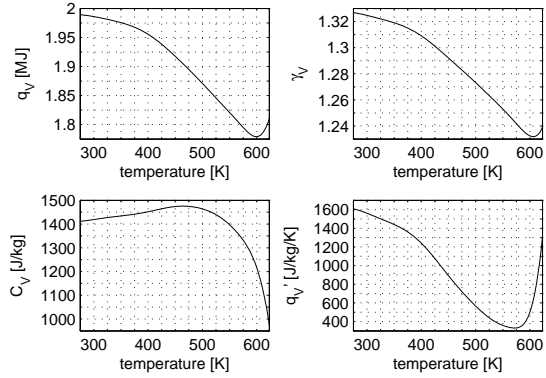}
\caption{{\footnotesize{Local optimal parameters for the vapor phase}}}
\label{fig:vap_para}
\end{figure}
The results are given in Figure \ref{fig:vap_para}.
\subsubsection{Liquid phase}\label{parl}
For the liquid phase we disclaim the simplifying assumption for the parameter $\pi$. Accordingly we 
are looking for 5 parameters. Therefore beside the relations 
(\ref{e:eos-e}) - (\ref{e:eos-s}) we can use a further relation. From thermodynamics it is known, that 
\[
 T\frac{\partial s(T,p)}{\partial T}=C^p
\]
with the {\em specific heat capacity at constant pressure}.
Using (\ref{e:eos-s}), (\ref{e:eos-T}) and (\ref{e:eos-e}) we find that
\begin{eqnarray*}
 C_L^p=T_L\frac{\partial s_L(T_L,p_L)}{\partial T_L}&=&C_L\gamma_L\\
 &=&\frac{1}{T_L}\left(\frac{p_L+\gamma_L\pi_L}{\rho_L(\gamma_L-1)}+\frac{\gamma_Lp_L-p_L}{\rho_L(\gamma_L-1)}\right)\\
 &=&\frac{1}{T_L}\left(e_L-q_L+\frac{p_L}{\rho_L}\right)\,.
\end{eqnarray*}
This gives us the further relation
\begin{equation}\label{e:eos-cp}
e_L=C^p_L T-\frac{p_L}{\rho_L}+q_L\,,
\end{equation}
where $C^p_L$ denotes the specific heat capacity at constant pressure of the liquid phase. 
Again we use the real equation of state 
to obtain the liquid density $\rho_{L*}$, 
the speed of sound $a_{L*}$, the entropy $s_{L*}$, the internal energy $e_{L*}$ as well as the specific heat capacity 
$C^p_{L*}$. Then from (\ref{e:eos-cp}) we find $q_L$. After that we calculate $\pi_L$ and $\gamma_L$ from (\ref{e:eos-e}) 
and (\ref{e:eos-a}). Finally we obtain $C_L$ from equation (\ref{e:eos-T}) and $q_L'$ from equation (\ref{e:eos-s}).
\begin{figure}[h!]
\includegraphics[scale=0.6]{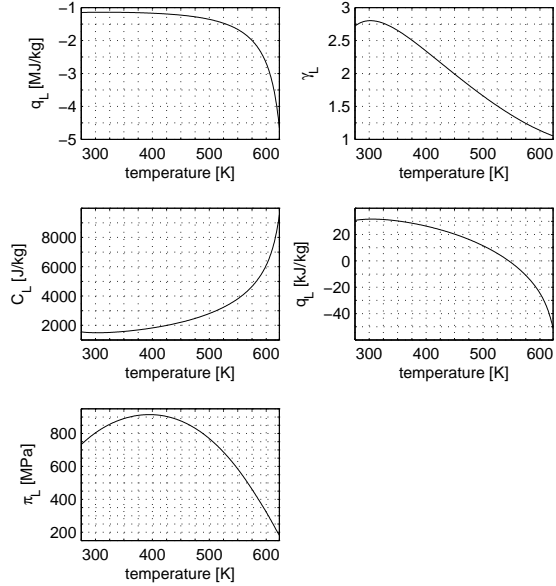}
\caption{{\footnotesize{Local optimal parameters for the liquid phase}}}
\label{fig:liq_para}
\end{figure}
The results are given in Figure \ref{fig:liq_para}.
\\\\
Using the parameters $\pi_V,\pi_L, \gamma_V,\gamma_L, C_V,C_L, q_V, q_L, q_V', q_L'$ obtained in Sections 
\ref{parv} and \ref{parl} the equations of state (\ref{e:eos-e})-(\ref{e:eos-s}) give the exact values for the densities,
the internal energies, the entropies and the sound speeds at saturation state $(p_*,T_*)$. 
Also we obtain the exact values for the Gibbs free energies and the enthalpies. From equations (\ref{e:gibbsrel})$_1$ and 
(\ref{e:gibbsrel})$_2$ as well as equation (\ref{e:f}) and the implicit function theorem we see, that we also find 
the exact value
\begin{equation}
T_{sat}'(p_*;T_*)\,.
\end{equation}
Due to the smoothness of all expressions for any given tolerance $\varepsilon>0$ we find a 
sufficiently small neighborhood of the 
saturation state $(p_*,T_*)$ such that all relevant physical states $\rho_k,e_k,s_k, a_k, g_k$ are approximated 
with a deviation less then $\varepsilon$.

\subsection{Proof of the statement for the real equation of state}
Assume, for any initial state in the vapor region, that is close to the saturation line, the corresponding wave curve 
in the $(p.T)$-phase space
and the (real) saturation line have the intersection point $(p_*,T_*)$. Assume further, that we used the optimal parameters 
$\pi_V,\pi_L, \gamma_V,\gamma_L, C_V,C_L, q_V, q_L, q_V', q_L'$, such that the equations of state 
(\ref{e:eos-e})-(\ref{e:eos-s}) give the exact values for $\rho_k,e_k,s_k, a_k, g_k$ for both phases 
at saturation state $(p_*,T_*)$. Then the following relation must hold at the intersection point
\begin{equation}\label{e:ablbed}
 \hat T'(p_*)\le T_{sat}'(p_*)\,.
\end{equation}
As before $\hat T(\hat p;p_*,T_*)$ denotes all admissible initial states $(\hat p,\hat T)$ in the vapor region, 
in a sufficiently small neighborhood of $(p_*,T_*)$. For the derivative we have
\begin{equation}\label{e:abl5}
  \hat T'(p_*;p_*,T_*)=\frac{T_*}{p_*}\cdot \frac{\gamma_V-1}{\gamma_V}\,.
\end{equation}
Moreover, we have
\begin{equation}\label{e:abl6}
 T_{sat}'(p_*;T_*)=T_*\cdot \frac{\frac{C_V(\gamma_V-1)}{p_*}-\frac{C_L(\gamma_L-1)}{p_*+\pi_L}}{C_V\gamma_V 
 -C_L\gamma_L+\frac{q_V-q_L}{T_*}}\,.
\end{equation}
In both equations (\ref{e:abl5}) and (\ref{e:abl6}) we used the local optimal parameters. As already explained at the 
end of Section \ref{approx}, equation (\ref{e:abl6}) gives the exact value for the derivative.
Simple estimations show, that 
\begin{equation}\label{e:abl7}
 \hat T'(p_*)> T_{sat}'(p_*)\,.
\end{equation}
Accordingly, there is no such intersection point. Due to the smoothness of all expressions and the exactness of 
(\ref{e:abl6}) this statement is true for the real equation of state and for all sufficiently good approximations 
of the real equation of state. If there is any set of parameters 
$\pi_V,\pi_L, \gamma_V,\gamma_L, C_V,C_L, q_V, q_L, q_V', q_L'$ such that (\ref{e:abl7}) is not satisfied, 
then the parameters obviously give a coarse approximation of the saturation line and the result is not meaningful. 
The same is true 
for any other choice of equations of state for the liquid and the vapor phase. We summarize
\begin{theorem}\label{n_allg}
 Using the real equations of state \cite{Wagner-Kruse:1998} or any good approximation of the real equation of state
 nucleation by compression cannot occur.
\end{theorem}

\section{Cavitation by expansion}\label{s:cav}
\begin{figure}[h!]
\includegraphics[scale=0.4]{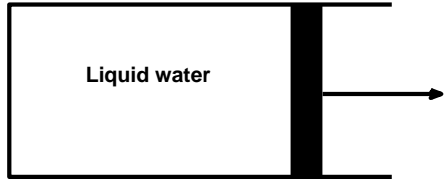}
\caption{{\footnotesize{Expansion of liquid water}}}
\label{expkolb}
\end{figure}

\begin{figure}[h!]
\includegraphics[scale=0.5]{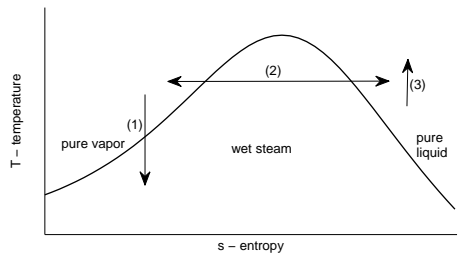}
\caption{{\footnotesize{Entropy temperature diagram: (1) isentropic expansion - rarefaction wave, (2) isothermal 
path, (3) isentropic compression}, see Iben \cite{I2}}}
\label{sdia}
\end{figure}

After the discussion of condensation by compression one may ask for the opposite case of cavitation by expansion. We will see, 
that this process is more complicated and we will distinguish between two cases. The liquid phase will be expanded 
in a manner 
that phase transition will occur. This case corresponds to a cavitation tube, which is a tube filled with liquid water and a 
flexible piston, see Figure \ref{expkolb}.
To illustrate the physics we give the $s-T$-diagram in Figure \ref{sdia}, where the path (1) corresponds to the process considered.
We have seen, that the Baer Nunziato type relaxation model allows the coexistence of both the vapor and the liquid phase 
in the same point of the physical domain at 
the saturation state. The same is true for the Eulerian model used by Dumbser et 
al.\ \cite{Dumbser}. Here the phase transition is modeled by an equilibrium assumption. The mixture is 
called {\em wet steam}. An expansion process such that a mixture of water vapor and liquid water (wet steam) is created, 
we call the process {\bf weak cavitation}. On the other hand, if pure water vapor is created, we call 
this process {\bf strong cavitation}.
\subsection{Cavitation in the weak sense}\label{weakcav}
If the liquid phase is expanded, a rarefaction wave will propagate through the liquid phase. In analogy to the condensation 
case we start with an arbitrary set of initial data in the liquid phase (Region 1). We construct the wave curve, that 
connects the initial state to all possible states behind the rarefaction wave. If (weak) cavitation, by sufficiently 
strong expansion, can occur this wave curve must have an intersection point with the saturation line. It is not 
surprising, that this is usually the case. Numerous examples can be found in the literature, see for instance the 
example {\em cavitation by strong rarefaction} using one set of Euler equations and the equilibrium assumption 
in Dumbser 
et al. \cite{Dumbser}. See further the {\em expansion tube problem} in Zein et al.\ \cite{Zein-Hantke-Warnecke:2010} using the Baer Nunziato 
type relaxation model. 
\subsection{Cavitation in the strong sense}\label{strongcav}
For the moment we restrict ourselves to the two model types, that allow the coexistence of vapor and liquid, 
\cite{Zein-Hantke-Warnecke:2010}, \cite{Dumbser}. Models, using a kinetic relation will be discussed later.

\begin{figure}[h!]
\includegraphics[scale=0.6]{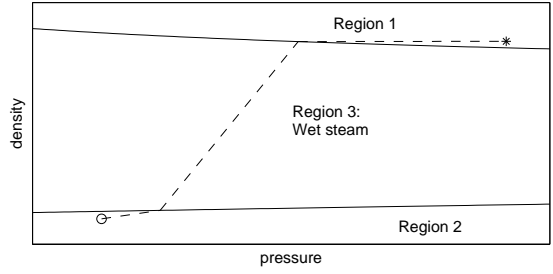}
\caption{{\footnotesize{Solid line: saturation curves $\rho_{V}(p_{sat})$, $\rho_{L}(p_{sat})$, 
dashed line: wave curve, star: initial state, circle: pure water vapor state}}}
\label{fig:cavwc}
\end{figure}

One may assume that for sufficiently strong expansions one may create pure water vapor. To illustrate that situation we 
refer to Figure \ref{fig:cavwc}.
Because Region 3, the wet steam region, reduces to a single line in the $(p,T)$ phase plan, we now plot 
all the data in the more descriptive $(p,\rho)$ phase plane. The wet steam region is bounded by the saturation (solid) 
lines. Assume, there is any initial liquid state, marked by the star and assume further, the liquid is expanded 
in a manner that pure water vapor is created. This state is indicated by the circle. 
Then there is a rarefaction wave curve, connecting the 
star and the circle state. This wave curve crosses the two curves $\rho_k(p_{sat})$. 
\begin{figure}[h!]
\includegraphics[scale=0.6]{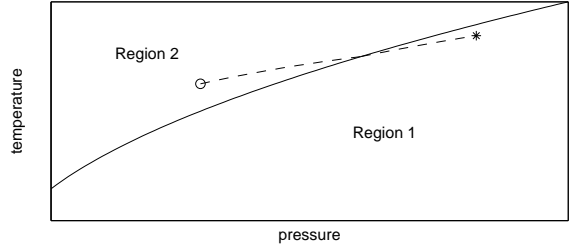}
\caption{{\footnotesize{Solid line: saturation curve, 
dashed line: wave curve, star: initial state, circle: pure water vapor state}}}
\label{fig:cav_sk}
\end{figure}
The same situation in the 
$(p,T)$ phase plane is given in Figure \ref{fig:cav_sk}. In the following we will show, that strong cavitation 
for the models considered cannot occur. The argument is similar to the argument in the condensation case. 
Assume, the circle state exists. This implies the existence of a rarefaction wave curve that connects both the 
circle and the star state. The intersection point of the wave curve 
and the saturation curve $T_{sat}(p)$ in the $(p,T)$ phase plane is called $(p_*,T_*)$. 
Let us consider that part of the wave curve, 
that is located in the vapor in Region 2. We denote this curve by $\hat T(p;p_*,T_*)$. As before 
for the intersection point we must have 
\[
 \hat T'(p_*;p_*,T_*)\le T'_{sat}(p_*)\,.
\]

The wave curve is found to be
\begin{equation}\label{wcrf}
 \hat T(p;p_*,T_*)=T_*\left(\frac{p}{p_*}\right)^{\frac{\gamma_V-1}{\gamma_V}}\,.
\end{equation}
This directly follows from Equation (\ref{e:eos-s}) and the fact, 
that the entropy is constant across a rarefaction wave.
Equation (\ref{wcrf}) implies 
\begin{equation}\label{ablrf}
 \hat T(p_*;p_*,T_*)=\frac{T_*}{p_*}\cdot \frac{\gamma_V-1}{\gamma_V}\,.
\end{equation}
This gives a contradiction, see (\ref{e:abl5}) and (\ref{e:abl7}).

\begin{theorem}
 Using the real equation of state or any good approximation of the real equation of state 
 strong cavitation by expansion cannot occur in an approach based on Euler equations and an equilibrium assumption, 
 \cite{Dumbser}, \cite{Zein-Hantke-Warnecke:2010}.
\end{theorem}

There is an alternative, very simple argument to show, that a pure liquid state and a pure water vapor state 
cannot be connected by a rarefaction wave. As already mentioned it is a well known fact, that the entropy is constant 
across a rarefaction wave. But it is also known, that for any temperature between triple point temperature 
and critical temperature $T_{tripel}=273.16K<T<T_{crit}=647.096K$ the entropies satisfy the inequalities
\begin{equation}\label{entropy}
 s_L < s_{crit} \qquad \mbox{and} \qquad s_{crit} < s_V\,.
\end{equation}
Obviously there are no liquid and vapor states with the same entropy. This also can be seen in Figure \ref{sdia}.

Using an Eulerian approach like Dumbser et al.\ \cite{Dumbser} one can obtain only weak cavitation. 
Wet steam is created, which is a mixture of water vapor and liquid water at saturation state. For the entropy 
one has
\begin{equation}\label{entrD}
 s_{mix}=\mu s_V+(1-\mu)s_L
\end{equation}
where $\mu\in[0,1]$ denotes the vapor mass fraction, see \cite{Dumbser}. Due to Equation (\ref{entropy}) 
the value $\mu$ is bounded for cavitation starting from pure liquid water. Using the steam tables of Wagner 
\cite{steam} we find
\begin{equation}\label{muest}
 \mu\le 0.5\,.
\end{equation}

In contrast to the equilibrium models discussed in the first part of Section \ref{strongcav} models 
using a kinetic relation are able to produce strong cavitation. This is clear by Section \ref{weakcav}.
A rarefaction wave curve in the liquid Region 1 can have an intersection point with the saturation line.
Here a critical state is reached. A cavitation criterion can be used, see Hantke et al.\ \cite{Hantke}. The liquid 
phase starts to evaporate and a pure vapor phase is created. An important difference to the previous models 
is, that the solution is nonsmooth and entropy production by phase transition is allowed. Therefore there is no contradiction to previous results.


\section{Conclusions}\label{s:clos}
We have seen, that condensation by compression cannot occur in an Eulerian approach. 
Due to the compression of the water vapor not only the pressure but also the temperature is increasing. 
Because of the temperature rise the saturation pressure is increasing. 
The key point is, that the saturation pressure increases much faster 
than the pressure inside the vapor phase. Therefore phase transition cannot take place.
From observations of nearly adiabatic flows, see Figure \ref{expkolb}, \ref{comkolb} 
this phenomenon is known.

\begin{itemize}
\item
The Euler equations correctly reflect this behavior. We gave a mathematical proof for this.
\item
Adiabatic compression of vapor does not lead to liquid phase, see Figure \ref{sdia}, path (1). To reach
this state, negative heat flow has to be used. This is equivalent to the use of
isothermal Euler equations such as done by Hantke et.\ al \cite{Hantke}. This corresponds to path (1) in Figure \ref{sdia}.
\end{itemize}

This effect comes up in much weaker form in the case of cavitation. The reason may be that the temperature 
changes due to the expansion are much smaller.

\begin{itemize}
\item
It is not possible to evaporate a pure liquid by a rarefaction wave completely, only with external supply of energy 
see again Figure \ref{sdia}. This is also equivalent to
the use of isothermal Euler equations \cite{Hantke}, see Figure \ref{sdia}, path (2).
\item
We gave a mathematical proof for this.
\end{itemize}

On the other hand this shows, that heat flow plays an important role in cavitating processes or in condensation 
processes caused by compression. Eulerian models are not appropriate to describe such effects.

\bibliographystyle{amsplain}

\end{document}